\documentclass[10pt,reqno]{amsart}
\usepackage{geometry} 
\usepackage{graphicx}
\usepackage{amsmath}
\usepackage{color}
\newtheorem{thm}[equation]{Theorem}
\newtheorem{pro}[equation]{Proposition}
\newtheorem{cor}[equation]{Corollary}
\newtheorem{lem}[equation]{Lemma}

\theoremstyle{definition}

\newtheorem{DEF}[equation]{Definition}
\newtheorem{rem}[equation]{Remark}

\usepackage{amsfonts,amssymb,amscd,amsmath,enumerate,verbatim,calc}
\def\span{\hbox{span}}
\def\aa{\mathcal A}

\def\bbbf{\mathbb{F}}

\def\LL{\mathcal{L}}

\def\sub{\subseteq}

\def\Lam{\Lambda}

\def\qed{\hfill$\Box$\vspace{3mm}}

\def\ad{\hbox{ad}}

\def\bbbz{{\mathbb Z}}
\def\bbbn{{\mathbb N}}

\def\vv{{\mathcal V}}
\def\uu{{\mathcal U}}

\def\Lam{\Lambda}
\def\uu{{\mathcal U}}
\begin{document}

\markboth{On the structure of graded $3-$Leibniz algebras} {V. Khalili}

%{E-mail: ......}
\date{}

\centerline{\bf On the structure of graded $3-$Leibniz algebras}

\vspace{.5cm}\centerline{Valiollah Khalili\footnote[1]{Department
of mathematics, Faculty of sciences, Arak University, Arak 385156-8-8349, Po.Box: 879, Iran. 
 V-Khalili@araku.ac.ir\\
\hphantom{ppp}2000 Mathematics Subject Classification(s): 17A32, 17A60, 17B22, 17B65.
\\
\hphantom{ppp}Keywords: graded Leibniz algebra, $3-$Leibniz algebra, structure theory }\;\;}

\vspace{1cm} \noindent ABSTRACT: We study the structure of a $3-$Leibniz algebra $T$ graded by an arbitrary
abelian group $G,$ which is considered of arbitrary dimension and over an arbitrary base field $\bbbf.$ We show that $T$ is of the form $T=\uu\oplus\sum_jI_j,$  with $\uu$ a linear subspace of $T_1,$ the homogeneous component  associated to the unit element $1$ in $G,$ and any  $I_j$ a well described graded
ideal of $T,$ satisfying 
$$
[I_j, T, I_k] = [I_j, I_k, T] = [T, I_j, I_k] = 0,
$$
if $j\neq k.$ In the case of $T$ being of maximal length, we
characterize the gr-simplicity of the algebra in terms of connections in the support of the
grading.  

\vspace{1cm} \setcounter{section}{0}
\section{Introduction}\label{introduction}

In \cite{L}, Loday introduces Leibniz algebras which are a “non-antisymmetric” generalization of Lie algebras. They play an important role in Nambu mechanics \cite{DT, G1}. Algebras endowed with an $n-$ary operation play important roles, among others, in Lie and Jordan theories, geometry, analysis and physics. Leibniz $n-$algebras
and its corresponding skew-symmetric version, named as Lie $n-$algebras or Filippov
algebras, arose in the setting of Nambu mechanics, a generalization of the Hamiltonian mechanics \cite{DT, N, T}.  Motivated by the ideas of Leibniz algebras,  the $n-$ary version of these algebras introduced \cite{CLP}, which appeared recently in the Bagger-Lambert-Gustavsson theory of multiple M2-branes \cite{BL, G2, HM}. For $n=2$ one recovers Leibniz algebras. Any Leibniz algebra  is also a Leibniz $n-$algebra under the following n-bracket:
$[x_1, x_2, ... , x_n]= [x_1, [x_2,... , [x_{n-1}, x_n] ...]].$
Conversely,  a Leibniz $(n + 1)-$algebra, can be made into a Leibniz algebra. This construction goes back to Gautheron \cite{G1} and plays an important role in \cite{DT}.
The particular case $n=3$ has found applications in string theory and M-branes \cite{BL, PS} and in the M-theory generalization of the Nahm’s equation proposed by Basu and Harvey \cite{BH}. It can also be used to construct solutions of the Yang-Baxter equation \cite{O}, which first appeared in statistical mechanics \cite{B}.

On the other hand, the descriptions of all possible gradings on algebras plays an important role both in the structure theory of finite dimentional and infinite dimentional  algebras and its applications (see \cite{BSZ, BT, C2, CC1, CC3, CF, CO, CS1, CS2, Kh, PZ}). In particular, the interest in group grading on $n-$ary algebras has been remarkable in the last years, motivated in part by their application in physics, geometry and topology where they appear as the natural framwork for an algebraic model (see \cite{BCNS, C3, CC2}). 

In the present paper we are going to the study of arbitrary graded $3-$Leibniz algebras
(not necessarily simple or finite-dimensional) introduced as the natural extension of
 graded Lie algebras and graded Leibniz algebras over an arbitrary base field $\bbbf$ by focussing on their structure (see \cite{C2, CD}).
We extend the techniques of connections in the support of $G$ developed for
graded Leibnize algebras in \cite{CD} to the framework of graded $3-$Leibniz algebras. It is worth mentioning that the so-called  connection techniques were first introduced by Calderon Martin in \cite{C1} for study of the structure of split Lie algebras with symmetric root systems.

This paper proceeds as follows; In section 2 we recall the main definitions and results related to the graded $3-$Leibniz algebras theory.  In Section 3 we develop the techniques of connections  in this framework, which become the main tool
in our study of graded $3-$Leibniz algebras.  We will show that an arbitrary $3-$Leibniz algebra $T$ is of the form $T=\uu\oplus\sum_jI_j,$  with $\uu$ a linear subspace of $T_1,$ the homogeneous component  associated to the unit element $1$ in $G,$  and any  $I_j$ a well described graded
ideal of $T,$ satisfying 
$$
[I_j, T, I_k] = [I_j, I_k, T] = [T, I_j, I_k] = 0,
$$
if $j\neq k.$ Finally, in section 4  we focus on the gr-simplicity and  gr-primeness of graded $3-$Leibniz algebras by centering
our attention in those of maximal length. Our results extend the ones
for graded Leibniz algebras \cite{CD}, graded Leibniz triple systems \cite{CC2} and for split $3-$Leibniz algebras \cite{CO}. Finally,
Section 5 provides a concrete example which characterizes the inner structure of
 graded $3-$Leibniz algebras.

Throughout this paper, algebras and vector spaces are over a field $\bbbf$ of characterestic zero,  We also consider an arbitrary $3-$Leibniz 
algebra (not necessarily simple or finite dimensional) over
 field $\bbbf$  and  an arbitrary abelian group $G.$

\section{Preliminaries}\setcounter{equation}{0}\

We begin by recalling the necessary background  on graded $3-$Leibniz algebra.

\begin{DEF}(see \cite{CLP}.) A  left Leibniz algebra $\L$ is a vector space over a field $\bbbf$ endowed with a bilinear product
$[., .]$ satisfying the Leibniz identity
$$
[x, [y, z]] = [[x, y], z]+[y, [x, z]],
$$
for all $x, y, z \in \L.$
\end{DEF}

\begin{DEF} (see \cite{CLP}.)  A  $3-$Leibniz algebra is a vector space $T$ , over an arbitrary base field $\bbbf$, with a trilinear operation $[., ., .]~:~ T\times T\times  T \longrightarrow T,$ called the triple product of $T,$ satisfying the following identity:
\begin{equation}\label{100}
[x, y, [a, b, c]] = [[x, y, a], b, c] + [a, [x, y, b], c] + [a, b, [x, y, c]],
\end{equation}
for any $x, y, a, b, c \in T.$ 
\end{DEF}

If $L$ is a Leibniz algebra with product $[., .],$ then $L$ becomes a $3-$Leibniz algebra by putting $[x, y, z] := [[x, y], z]$.  $3-$Lie or a $3-$Nambu-Lie algebra is a $3-$Leibniz algebra such that trilinear operation is skew-symmetric. Another big class of $3-$Leibniz algebras  are the so called Lie triple systems (see \cite{CLP}).
 
Let $T$ be a $3-$Leibniz algebra. A linear subspace of $T$ closed under the triple product is called a $3-$Leibniz {\em subalgebra} of $T.$ 
An {\em ideal} of $T$ is a linear subspace $I$ which satisfies $[I, T , T ] + [T , I, T ] + [T , T , I] \subset I.$ The {\em annihilator} of $T$ is the set
$$
Ann(T) = \{x \in T~:~  [x, T , T ] + [T , x, T ] + [T , T , x] = 0\}.
$$
It is straightforward to check that $Ann(T)$ is an ideal of $T$ (see \cite{CO}).

The ideal $J$ generated by the set
$$
\{[x, x, y], [x, y, x], [y, x, x]~:~ x, y \in T\},
$$
plays an important role in the theory since it determines the (possible) non-Lie character of $T$. Clearly, $T$ is a $3-$Lie algebra if and only if $J =\{0\}.$  From the Eq. (\ref{100}), this ideal also satisfies
\begin{equation}\label{101}
[T, T,  J]=[T,J,T] = 0.
\end{equation}

We will say that $T$ is a simple $3-$Leibniz algebra if $[T, T, T]\neq 0$ and its only ideals are $\{0\}, J$ and $T$ (see \cite{CC2}). 

Recall that, for a $3-$Leibniz algebra $T,$  a linear map $D~:~T\longrightarrow T$ is said to be a ternary derivation if it satisfies
$$
D([x, y, z]) = [D(x), y, z] + [x, D(y), z] + [x, y, D(z)],
$$
for all $x, y, z\in T.$

Using this terminology, the identity (\ref{100}) says that the left multiplication operators $\ad(x, y)~:~ T\longrightarrow T$ given by $\ad(x, y)(z) :=
[x, y, z]$ are ternary derivations.
The identity (\ref{100}) applies to get that $\LL := span_{\bbbf}\{ad(x, y)~:~ x, y \in T \}$ with the bracket
$$
[\ad(a, b), \ad(c, d)]= \ad([a, b, c], d)+ ad(c, [a, b, d])~~~\forall a, b, c, d\in T,
$$
is a Lie algebra (see \cite{CO}). 

\begin{DEF}\cite{CO} The standard embedding of a $3-$Leibniz algebra $T$ is the $2-$graded algebra $\aa = \aa^0\oplus\aa^1,$  where $\aa^0=\LL,~\aa^1=T,$ and product given by
\begin{eqnarray*}
(\ad(x, y), z).(\ad(u, v), w)&:=& (\ad([x, y, u], v) + \ad(u, [x, y, v]) + \ad(z, w)\\
&,&[x, y, w]-[u, v, z]).
\end{eqnarray*}
\end{DEF}

Let us observe that $\aa^0$ with the product induced by the one in $\aa = \aa^0\oplus\aa^1$ becomes a Lie algebra, $\aa$ is not, in general, a ($2-$graded) Lie algebra.

Next, we introduce the class of graded algebras in the framework of $3-$Leibniz algebras. Observe that the product in $\aa$ gives us a natural action:
$$
\aa^0\times\aa^1\longrightarrow\aa,~~~(x, y)\longmapsto xy.
$$ 
Since $\aa^0$ is a Lie algebra, the identity (\ref{100}) allows us to conclude that the action above endows $\aa^1$ with an $\aa^0-$module structure. Thus, we can introduce the concept of graded $3-$Leibniz algebra (graded $3-$Lie algebras which not yet charatrized) in a similar spirit to the ones of  graded Lie triple system \cite{C3}.

\begin{DEF} Let $T$ be a $3-$Leibniz algebra. It is said that $T$ is graded by means
of an abelian group $G$ if it decomposes as the direct sum of linear subspaces
$$
T=\bigoplus_{g\in G} T_g,
$$
where the homogeneous components satisfy $[T_g, T_h, T_k] \subset T_{ghk}$ for any $g, h, k \in	G$ (denoting by juxtaposition the product and unit element 1 in $G$). We call the $G-$support of the grading the set
$$
\Sigma^1~:=\{g\in G\setminus\{1\}~:~T_g\neq 0\}.
$$
\end{DEF}

We will also say that the $G-$support of the grading is symmetric if $g\in\Sigma^1$ implies
$g^{-1}\in\Sigma^1.$ We finally note that graded $3-$Lie algebras, split $3-$Leibniz algebras and graded Leibniz triple systems are examples of graded $3-$Leibniz algebras and so the present paper extends the results in \cite {CC2, CO}.

As it is usual in the theory of graded algebras, the  regularity conditions will be understood in the graded sense compatible with the $3-$Leibniz algebra structure. That
is, a $3-$Leibniz  subalgebra (or ideal) of $T$ is a graded linear subspace $S$ (or $I$).   Moreprecisely, $S$ (or $I$) splits as $S=\bigoplus_{g\in G} S_g,~~S_g=S\cap T_g,$ similarly for $I.$ Also we will say that $T$ is a gr-simple $3-$Leibniz algebra if $[T, T, T]\neq 0$ and its only graded ideals are $\{0\}, J$ and $T.$  

Recall that, let $L$ be an arbitrary Lie algebra over $\bbbf.$ As usual, the term grading will always mean abelian group grading, that is, a decomposition in linear subspaces $L=\bigoplus_{g\in G} L_g,$ where $G$ is an abelian group and the homogeneous spaces satisfy $[L_g, L_h] \subset L_{gh}.$ We also call the support of the grading the set 
$$
\{g\in G\setminus\{1\}~:~L_g\neq 0\}.
$$

\begin{pro}\label{.} Let $T$ be a graded $3-$Leibniz algebra and let $\aa = \LL\oplus T$
be its standard embedding algebra. Then $\LL$ is a $G-$graded Lie algebra.
\end{pro}
\noindent {\bf Proof.} The proof is vertically identical to the
proof  of Proposition 2.1 in \cite{C3}. \qed

Observe that for any $g, h \in G$ we have
$$
[\LL_g, \LL_h] \subset\LL_{gh}.
$$
In the following, we shall denote by $\Sigma^0$ 
the support of the $G-$graded Lie algebra $\LL.$

\section{$\Sigma^1-$Connections and Decompositions} \setcounter{equation}{0}\

In this section,  we are going to develop connection in the support
techniques for graded  $3-$Leibniz algebras as the main tool in our study. From now on and throughout the paper, $T$ denotes an arbitrary   $G-$graded $3-$Leibniz algebra with the  standard embedding  $\aa = \LL\oplus T,$  a
$G-$support $\Sigma^1,$ and
\begin{equation}\label{0'}
T=\bigoplus_{g\in G} T_g=T_1\oplus(\bigoplus_{g\in\Sigma^1} T_g),
\end{equation}
the corresponding grading. Denote by $-\Sigma^i=\{-g~:~g\in\Sigma^i\}, i=0, 1.$

\begin{DEF}\label{conn}
Let $g, h\in\Sigma^1.$ We sall say that
{\em $g$ is connected to $h$} and denote it by $g\sim h$ such that
if there exists a family
$$
\{g_1, g_2, g_3, ...,
g_{2n+1}\}\subset\pm\Sigma^1\cup\{1\},
$$
satisfying the following conditions;\\

\begin{itemize}
\item[(1)] $\{g_1, g_1g_2g_3, ..., g_1g_2g_3 ... g_{2n}g_{2n+1}\}\subset\pm\Sigma^1,$\\

\item[(2)] $\{g_1g_2, g_1g_2g_3g_4, ... , g_1g_2g_3 ... g_{2n}\}\subset\pm\Sigma^0,$\\

\item[(3)] $g_1 = g$ and $g_1g_2g_3 ... g_{2n}g_{2n+1} \in\{h, h^{-1}\}.$
\end{itemize}
The family $\{g_1, g_2, g_3, ..., g_{2n+1}\}$ is
called a {\em $\Sigma^1-$connection} from $g$ to $h.$
\end{DEF}

The next result shows that the connection relation in $\Sigma^1$ is an equivalence relation. 

\begin{pro}\label{rel}
The relation $\sim$ in $\Sigma^1,$ defined by
$g\sim h$ if and only if $g$ is connected to $h,$
is an equivalence relation.
\end{pro}
\noindent {\bf Proof.} The proof is analogous to the one for graded Lie triple systems given in \cite{C3}, Proposition 3.1.  \qed

By the above proposition we can consider the quotient set
$$
\Sigma^1/\sim :=\{[g] :
g\in\Sigma^1\},
$$
then $[g]$ is the set of  elements  in the $G-$support of the grading which are connected to $g.$ By the definition of $\sim,$ it is clear that if $h \in[g]$ and $h^{-1}\in\Sigma^1,$ then $h^{-1}\in[g].$

Our  goal in this section is to associate an ideal
$I_{[g]}$ of the graded $3-$Leibniz algebra  $T$ to any $[g]\in\Sigma^1/\sim.$ Fix
$g\in\Sigma^1,$   we  define
\begin{equation}\label{0.0}
I_{1, [g]} :=\span_{\bbbf}\{[T_h, T_k, T_{(hk)^{-1}} ] ~: h\in[g], k\in [g]\cup\{1\} \}\sub T_1,
\end{equation}
and
\begin{equation}\label{0.0.}
\vv_{[g]} :=\bigoplus_{h\in[g]}T_h.
\end{equation}
Finally, we denote by $I_{[g]}$ the direct sum of  two 
subspaces above, that is,
\begin{equation}\label{0.00}
I_{[g]} :=I_{1, [g]}\oplus\vv_{[g]}.
\end{equation}

\begin{pro}\label{subalg} For any
$g\in\Sigma^1,$ the graded subspace $I_{[g]}$ is a graded
subalgebra of $T.$ We will refer to $I_{[g]}$ as the $3$-Leibniz subalgebra
of $T$ associated to $[g].$
\end{pro}
\noindent {\bf Proof.} We are going to check that
$I_{[g]}$ satisfies $[I_{[g]}, I_{[g]}, I_{[g]}]\subset I_{[g]}.$ We have
\begin{eqnarray}\label{200}
\nonumber[I_{[g]}, I_{[g]}, I_{[g]}]&=&[I_{1, [g]}\oplus\vv_{[g]}, I_{1, [g]}\oplus\vv_{[g]}, I_{1, [g]}\oplus\vv_{[g]}]\\
&\subset&[I_{1, [g]}, I_{1, [g]}, I_{1, [g]}]+[I_{1, [g]}, I_{1, [g]},  \vv_{[g]}]\\
\nonumber&+&[I_{1, [g]}, \vv_{[g]}, I_{1, [g]}]+[\vv_{[g]}, I_{1, [g]}, I_{1, [g]}]+[I_{1, [g]}, \vv_{[g]}, \vv_{[g]}]\\
\nonumber&+&[\vv_{[g]}, I_{1, [g]}, \vv_{[g]}]+[\vv_{[g]}, \vv_{[g]}, I_{1, [g]}]
+[\vv_{[g]}, \vv_{[g]}, \vv_{[g]}].
\end{eqnarray}
Since $I_{1, [g]}\subset T_1$ we clearly have 
\begin{equation}\label{103}
[I_{1, [g]}, I_{1, [g]}, I_{1, [g]}]\subset[I_{1, [g]},, T_1, T_1]\subset I_{1, [g]} \subset I_{[g]}, 
\end{equation}
and also 
\begin{equation}\label{104}
[I_{1, [g]}, I_{1, [g]},  \vv_{[g]}]+[I_{1, [g]}, \vv_{[g]}, I_{1, [g]}]+[\vv_{[g]}, I_{1, [g]}, I_{1, [g]}]\subset\vv_{[g]}\subset I_{[g]}. 
\end{equation}
Let us consider the sammand $[I_{1, [g]}, \vv_{[g]}, \vv_{[g]}]$ in (\ref{200}). We have 
$$
[I_{1, [g]}, \vv_{[g]}, \vv_{[g]}]\subset [T_1, \vv_{[g]}, \vv_{[g]}].
$$ 
suppose $[T_1, T_h, T_k]\neq 0$ for some  $h, k\in [g],$  we get  $h\in\Sigma^0$ and $hk\in\Sigma^1\cup\{1\}.$ From here, if $hk = 1$ clearly $[T_1, T_h, T_{h^{-1}}]\subset I_{1,[g]}.$ If $hk\neq 1$  and $\{g_1, ..., g_{2n+1}\}$ is a connection from $g$ to $h$ then $\{g_1, ..., g_{2n+1}, 1, k\}$ is a connection from $g$ to $hk$ in case $g_1g_2 ... g_{2n+1} = h$ and $\{g1, ... , g_{2n+1}, 1, k^{-1}\}$ in case $g_1g_2 ... g_{2n+1} = h^{-1}$ being so $hk \in [g],$ hence $[T_1, T_h, T_k]\subset \vv_{[g]}.$ In any case, we have proved $hk\in [g]$ so $[T_1, T_h, T_k]\subset T_{hk}\subset \vv_{[g]}.$ Therefor,
\begin{equation}\label{106}
[I_{1, [g]}, \vv_{[g]}, \vv_{[g]}]\subset I_{[g]}.
\end{equation}
A similar argument as  the previous case, we obtain
\begin{equation}\label{109}
[\vv_{[g]}, I_{1, [g]}, \vv_{[g]}]\subset I_{[g]}.
\end{equation}
Next, consider the sammand $[\vv_{[g]}, \vv_{[g]}, I_{1, [g]}]$ in (\ref{200}). We have 
$$
[\vv_{[g]}, \vv_{[g]}, I_{1, [g]}]\subset [\vv_{[g]}, \vv_{[g]}, T_1].
$$ 
Assume now $[T_h, T_k, T_1]\neq 0,$ for some  $h, k\in [g]$ then $hk\in\Sigma^0$ and $hk\in\Sigma^1\cup\{1\}.$  If $hk = 1$ clearly $[T_h, T_k, T_1]\subset I_{1,[g]}.$ If  $hk\neq 1$ then  $hk\in\Sigma^0\cap\Sigma^1.$ Consider $\{g_1, ..., g_{2n+1}, k, 1\}$  a connection from $g$ to $h.$ We can also consider this is a connection from $g$ to $hk.$ More precisely,  $\{g_1, ..., g_{2n+1}, k, 1\}$  is a connection from $g$ to $hk$ in case $g_1g_2 ... g_{2n+1} = h$ and $\{g1, ... , g_{2n+1}, k^{-1}, 1\}$ in case $g_1g_2 ... g_{2n+1} = h^{-1},$  being so $hk \in [g],$ hence $[T_1, T_h, T_k]\subset \vv_{[g]}.$ In any case, we have proved $hk\in [g].$ Thus, we have $[T_1, T_h, T_k]\subset T_{hk}\subset \vv_{[g]}.$ Therefor,
\begin{equation}\label{113}
[\vv_{[g]}, \vv_{[g]}, I_{1, [g]}]\subset I_{[g]}.
\end{equation}
Finally, let us consider the last summand in  Eq. (\ref{200}) and suppose
there exist $h, k, l\in [g]$ such that $[T_h, T_k, T_l]\neq 0$ then  $hkl\in\Sigma^1\cup\{1\}.$  If  $hkl = 1$ then $[T_h,T_k, T_l]\subset T_1\subset \vv_{[g]}.$  Suppose now that $hkl\in\Sigma^1,$  and consider two
cases. In case $hk=1,$  we obtain that $[T_h,T_k, T_l]\subset T_l\subset \vv_{[g]}.$ If $hk\neq 1$ then we get that $hk\in\Sigma^0$ and $hkl\in\Sigma^1.$ If $\{g_1, ..., g_{2n+1}\}$ is a connection from $g$ to $h.$ We clearly have $\{g_1, ... , g_{2n+1}, k, l\}$ is a connection from $g$ to $hkl$ in case $g_1 ... g_{2n+1} = h$ and $\{g_1 ... g_{2n+1}, k^{-1}, l^{-1}\}$ it is in case $g_1 ... g_{2n+1} = h^{-1}.$ Thus we have get  $hkl\in [g]$  and so $[T_h, T_k, T_l]\subset \vv_{[g]}.$ Therefor,
\begin{equation}\label{117}
[\vv_{[g]}, \vv_{[g]}, \vv_{[g]}]\subset I_{[g]}.
\end{equation}
From Eqs. (\ref{103})-(\ref{117}) we conclud that $I_{[g]}$ is a 
subalgebra of graded $3-$Leibniz algebra $T.$\qed

In the following, we are going to show that  for any $g\in\Sigma^1,$ the graded subalgebra
of $T$ associated to $[g]$ is actually a graded ideal of $T.$ We need to state some preliminary results, which they can be
proved as in \cite{C3, CC2}.

\begin{lem}\label{105} The following assertions hold.
\begin{itemize}
\item[(1)] If $g, h \in\Sigma^1$ with $gh\in\pm\Sigma^0\cup\{1\}$ then $h\in [g].$

\item[(2)] If $g, h \in\Sigma^1$ and $g\in\pm\Sigma^0$ with $gh\in\pm\Sigma^1\cup\{1\}$ then $h\in [g].$

\item[(3)] If $g, h \in\Sigma^1$ such that $\bar{h}\notin [g],$ then $[T_g, T_{\bar{h}}]=[\LL_g, T_{\bar{h}}]=[\LL_g, \LL_{\bar{h}}]=0.$

\item[(4)] If $g, \bar{h} \in\Sigma^1$ are not connected, then $[T_g, T_{g^{-1}}, T_{\bar{h}}]=0.$
\end{itemize}
\end{lem}
\noindent {\bf Proof.} The proof is analogous to the one for graded Lie triple systems given in \cite{C3} Lemma 4.1 and 4.2.\qed

\begin{lem}\label{4.3} For any $g_0\in\Sigma^1,$ if $g\in [g_0]$ and $h, k\in\Sigma^1\cup\{1\},$ the following assertions hold.
\begin{itemize}
\item[(1)] If $[T_g, T_h, T_k]\neq 0$ then	$h, k, ghk\in [g_0]\cup\{1\}.$

\item[(2)] If $[T_h, T_k, T_g]\neq 0$ then	$h, k, ghk\in [g_0]\cup\{1\}.$

\item[(3)] If $[T_h, T_g, T_k]\neq 0$ then	$h, k, ghk\in [g_0]\cup\{1\}.$
\end{itemize}
\end{lem}
\noindent {\bf Proof.} See Lemma 4.3 in \cite{C3}.\qed

\begin{lem}\label{4.4} For any $g_0\in\Sigma^1,$ if $g, k\in [g_0]$ and $h\in\Sigma^0\cup\{1\}$ with  $ghk=1$ and $l, m\Sigma^1\cup\{1\},$ the following assertions hold.
\begin{itemize}
\item[(1)] If $[[T_g, T_h, T_k], T_l, T_m] = 0,$ then $l, m, lm \in [g_0] \cup\{1\}.$
\item[(2)] If $[T_l, [T_g, T_h, T_k], T_m] = 0,$ then $l, m, lm \in [g_0] \cup\{1\}.$
\item[(3)] If $[T_l, T_m, [T_g, T_h, T_k]] = 0,$ then $l, m, lm \in [g_0].$
\end{itemize}		
\end{lem}
\noindent {\bf Proof.} See Lemma 4.4 in \cite{C3}.\qed

\begin{lem} For any $g_0\in\Sigma^1,$ if $g, k\in [g_0], h\in\Sigma^0\cup\{1\}$ with  $ghk=1$ and $\bar{h}\in [g_0],$ the following assertions hold.
\begin{itemize}
\item[(1)]  $ [[T_g, T_h, T_k], T_h] = 0.$ 
\item[(2)] $[[T_g, T_h, T_k],\aa^0_{\bar{h}}] = 0.$
\item[(3)]  $[[T_g, T_h, T_k], T_1, T_{\bar{h}}] = 0.$
\end{itemize}		
\end{lem}
\noindent {\bf Proof.} See Lemma 4.5 in \cite{C3}.\qed

\begin{thm}\label{main1} The following assertions hold
\begin{itemize}
\item[(1)] For any $g_0 \in\Sigma^1$ the  subalgebra
$$
I_{[g_0]} =I_{1, [g_0]}\oplus\vv_{[g_0]},
$$
of $T$ associated to $[g_0]$ is an ideal of $T.$
\item[(2)] If $\LL$ is simple, then $\Sigma^1$ has all of its elements connected and
$$
T_1=\sum_{g\in\Sigma^1,~h\in\Sigma^1\cup\{1\}}[T_g, T_h, T_{(gh)^{-1}}].
$$
\end{itemize}
\end{thm}
\noindent {\bf Proof.}

(1) We are going to check 
\begin{equation}\label{107}
[I_{[g_0]}, T, T] + [T, T_{[g_0]}, T] + [T, T, I_{[g_0]}]\subset I_{[g_0]}.
\end{equation}
Let us consider the first sammand in the left hand in Eq. (\ref{107}). By Eqs. (\ref{0'}) and (\ref{0.00}), we have
\begin{eqnarray}\label{108}
[I_{[g_0]}, T, T]&=&[I_{1, [g_0]}\oplus\vv_{[g_0]},T_1\oplus(\bigoplus_{g\in\Sigma^1} T_g), T_1\oplus(\bigoplus_{h\in\Sigma^1} T_h)]\\
\nonumber&\subset& [I_{1, [g_0]}, T_1, T_1]+[I_{1, [g_0]}, T_1, \bigoplus_{h\in\Sigma^1} T_h]+[I_{1, [g_0]}, \bigoplus_{g\in\Sigma^1} T_g, T_1]\\
\nonumber&+&[I_{1, [g_0]}, \bigoplus_{g\in\Sigma^1} T_g, \bigoplus_{h\in\Sigma^1} T_h]+[\vv_{[g_0]}, T_1, T_1]+[\vv_{[g_0]}, T_1, \bigoplus_{h\in\Sigma^1} T_h]\\
\nonumber&+&[\vv_{[g_0]}, \bigoplus_{g\in\Sigma^1} T_g, T_1]+[\vv_{[g_0]}, [\bigoplus_{g\in\Sigma^1} T_g, \bigoplus_{h\in\Sigma^1} T_h].
\end{eqnarray}
Taking into account Eq. (\ref{0.0}) we have by Eq. (\ref{103}) that 
\begin{equation}\label{109'}
[I_{1, [g_0]}, T_1, T_1]\subset I_{[g_0]}.
\end{equation}
By Eq. (\ref{0.0}) and Lemma \ref{4.4}, we get
$$
[I_{1, [g_0]}, T_1, T_h]+[I_{1, [g_0]}, T_g, T_1]+[I_{1, [g_0]}, T_g, T_h]\subset I_{[g_0]}, \forall g, h\in\Sigma^1.
$$
Therefore,
\begin{equation}\label{110}
[I_{1, [g_0]}, T_1, \bigoplus_{h\in\Sigma^1} T_h]+[I_{1, [g_0]}, \bigoplus_{g\in\Sigma^1} T_g, I_1]+[I_{1, [g_0]}, \bigoplus_{g\in\Sigma^1} T_g, \bigoplus_{h\in\Sigma^1} T_h]\subset I_{[g_0]}.
\end{equation}
Now, taking into account Eq. (\ref{0.0.}) and Lemma \ref{4.3}, we get
$$
[T_k, T_1, T_1]+[T_k, T_1, T_h]+[T_k, T_g, T_1]\subset I_{[g_0]}, \forall k\in [g_0], g, h\in\Sigma^1.
$$
Therefore,
\begin{equation}\label{111}
[\vv_{[g_0]}, T_1, T_1]+[\vv_{[g_0]}, T_1, \bigoplus_{h\in\Sigma^1} T_h]+[\vv_{[g_0]}, \bigoplus_{g\in\Sigma^1} T_g, T_1]\subset I_{[g_0]}.
\end{equation}
Again by Eq. (\ref{0.0.}) and Lemma \ref{4.3}, we have
$$
[T_k, T_g, T_h]\subset I_{[g_0]}, \forall k\in [g_0], g, h\in\Sigma^1,
$$
and so
\begin{equation}\label{112}
[\vv_{[g_0]}, \bigoplus_{g\in\Sigma^1} T_g, \bigoplus_{h\in\Sigma^1} T_h]\subset I_{[g_0]}.
\end{equation}
Acording to Eq. (\ref{108}), from Eqs. (\ref{109'})-(\ref{112}) we conclude that
$$
[I_{[g_0]}, T, T]\subset I_{[g_0]}.
$$
By Lemmas \ref{4.3} and \ref{4.4}, a similar argument as above gives us $[T, T_{[g_0]}, T]\subset I_{[g_0]}$ and $[T, T, I_{[g_0]}]\subset I_{[g_0]}.$ Consequently,  $I_{[g_0]}$ is an ideal of $T.$

(2) The simplicity of $T$ implies $I_{[g]}\in\{J, T\}.$  If for any $[g]\in\Sigma^1/\sim $ is such that $I_{[g]}=T.$ Then $[g] =\Sigma^1.$ Hence, $\Sigma^1$ has all of its elements connected. Otherwise, if
$I_{[g]}=J$ for any $[g]\in\Sigma^1/\sim $ then $T_{[g]}=T_{[h]}$ for any $g, h\in\Sigma^1$ and so $[g] = \Sigma^1,$ we also conclude that $\Sigma^1$ has all its  elements connected. Observe that $T_1=\sum_{g\in\Sigma^1,~h\in\Sigma^1\cup\{1\}}[T_g, T_h, T_{(gh)^{-1}}]$ follows in any case.\qed

\begin{thm}\label{main2} For a vector space complement $\uu$  in $T_{1, \Sigma^1} :=\span_{\bbbf}\{[T_g, T_h, T_{(gh)^{-1}} ] ~: g\in\Sigma^1, h\in\Sigma^1\cup\{1\} \},$ we have
$$
T=\uu\oplus\sum_{[g]\in\Sigma^1/\sim}I_{[g]},
$$
where any $I_{[g]}$ is one of the  ideals of $T$
described in Theorem \ref{main1}-(1), satisfying $[I_{[g]}, T,
I_{[h]}]=[I_{[g]}, I_{[h]}, T]+[T, I_{[g]}, I_{[h]}]=0,$ whenever $[g]\neq [h].$
\end{thm}
\noindent {\bf Proof.} Each $I_{[g]}$ is well defined and by
Theorem \ref{main1}-(1), an ideal of $T.$ By Eq. (\ref{0'}), we have that
$$
T=T_1\oplus(\bigoplus_{g\in\Sigma^1}T_g)=(\uu\oplus T_{1, \Sigma^1})\oplus\sum_{[g]\in\Sigma^1/\sim}I_{[g]}.
$$
It follows that $\bigoplus_{g\in\Sigma^1}T_g=\bigoplus_{[g]\in\Sigma^1/\sim}\vv_{[g]}$ and $T_{1, \Sigma^1}=\sum_{[g]\in\Sigma^1/\sim} I_{1, [g]},$ which imply
$$
T=(\uu\oplus T_{1, \Sigma^1})\oplus\sum_{[g]\in\Sigma^1/\sim}I_{[g]}=\uu\oplus\sum_{[g]\in\Sigma^1/\sim}I_{[g]}.
$$
Next, it is sufficient to show that $[I_{[g]}, T, I_{[h]}]=0$ if $[g] \neq [h].$ By Eqs. (\ref{0'}) and (\ref{0.00}), we have
\begin{eqnarray}\label{115}
	[I_{[g]}, T, I_{[h]}]&=&[I_{1, [g]}\oplus\vv_{[g]},T_1\oplus(\bigoplus_{k\in\Sigma^1} T_k), I_{1, [h]}\oplus\vv_{[h]}]\\
	\nonumber&\subset& [I_{1, [g]}, T_1,  I_{1, [h]}]+[I_{1, [g]}, T_1, \vv_{[h]}]+[I_{1, [g]}, \bigoplus_{k\in\Sigma^1} T_k,  I_{1, [h]}]\\
	\nonumber&+&[I_{1, [g]}, \bigoplus_{k\in\Sigma^1} T_k, \vv_{[h]}]+[\vv_{[g]}, T_1, I_{1, [h]}]+[\vv_{[g]}, T_1, \vv_{[h]}]\\
	\nonumber&+&[\vv_{[g]}, \bigoplus_{k\in\Sigma^1} T_k,  I_{1, [h]}]+[\vv_{[g]}, \bigoplus_{k\in\Sigma^1} T_k,  \vv_{[h]}].
\end{eqnarray}
If  $[g] \neq [h],$ by Lemmas \ref{4.3} and \ref{4.4}, it is easy to see
\begin{eqnarray}\label{116}
0&=&[I_{1, [g]}, T_1, \vv_{[h]}]=[I_{1, [g]}, \bigoplus_{k\in\Sigma^1}T_k, \vv_{[h]}]=[\vv_{[g]}, T_1, I_{1, [h]}]\\
\nonumber&=&[\vv_{[g]}, T_1, \vv_{[h]}]=[\vv_{[g]}, \bigoplus_{k\in\Sigma^1} T_k,  I_{1, [h]}]=[\vv_{[g]}, \bigoplus_{k\in\Sigma^1} T_k,  \vv_{[h]}].
\end{eqnarray}
Now, let us consider the summand $[I_{1, [g]}, \bigoplus_{k\in\Sigma^1}T_k, I_{1, [h]}]$ in Eq. (\ref{115}). Given $[T_{h_1}, T_{h_2}, T_{(h_1 h_2)^{-1}}]\in I_{1, [h]}$ with $h_1\in [h], h_2\in [h]\cup\{1\}.$ By Eq. (\ref{100}) we have
\begin{eqnarray*}
[I_{1, [g]}, \bigoplus_{k\in\Sigma^1}T_k, [T_{h_1}, T_{h_2}, T_{(h_1 h_2)^{-1}}]]&\subset&
[[I_{1, [g]}, \bigoplus_{k\in\Sigma^1}T_k, T_{h_1}], T_{h_2}, T_{(h_1 h_2)^{-1}}]\\
&+&[T_{h_1}, [I_{1, [g]}, \bigoplus_{k\in\Sigma^1}T_k, T_{h_2}], T_{(h_1 h_2)^{-1}}]\\
&+&[T_{h_1}, T_{h_2}, [[I_{1, [g]}, \bigoplus_{k\in\Sigma^1}T_k, T_{(h_1 h_2)^{-1}}]]
\end{eqnarray*}

Given  $[T_{g_1}, T_{g_2}, T_{(g_1 g_2)^{-1}}]\in I_{1, [g]}$ with $g_1\in[g], g_2\in [g]\cup\{1\}.$ By the fact that $g\nsim h$ and Lemma \ref{4.4}, we get
\begin{equation*}
	\begin{split}
		&[I_{1, [g]}, \bigoplus_{k\in\Sigma^1}T_k, T_{h_1}]=[[T_{g_1}, T_{g_2}, T_{(g_1 g_2)^{-1}}], \bigoplus_{k\in\Sigma^1}T_k, T_{h_1}]=0,\\
		& [I_{1, [g]}, \bigoplus_{k\in\Sigma^1}T_k, T_{h_2}]=[[T_{g_1}, T_{g_2}, T_{(g_1 g_2)^{-1}}], \bigoplus_{k\in\Sigma^1}T_k, T_{h_2}]=0,\\
		&[I_{1, [g]}, \bigoplus_{k\in\Sigma^1}T_k, T_{(h_1 h_2)^{-1}}]=[[T_{g_1}, T_{g_2}, T_{(g_1 g_2)^{-1}}], \bigoplus_{k\in\Sigma^1}T_k, T_{(h_1 h_2)^{-1}}]=0.
	\end{split}
\end{equation*}
Therefore,
\begin{equation}\label{117}
[I_{1, [g]}, \bigoplus_{k\in\Sigma^1} T_k, I_{1, [h]}]=0.
\end{equation}
A similar argument gives us
\begin{equation}\label{118}
[I_{1, [g]}, T_1,  I_{1, [h]}]=0.
\end{equation}
From Eqs. (\ref{116})-(\ref{118}), we proved that $[I_{[g]}, T, I_{[h]}]=0$ if $[g] \neq [h].$ A similar discussion as above gives us $[I_{[g]}, I_{[h]}, T]=[T, I_{[g]}, I_{[h]}]=0,$ whenever $[g]\neq [h].$\qed 

\begin{cor}\label{2.17} If $Ann(T)=0$ and   $T_1=\span_{\bbbf}\{[T_g, T_h, T_{(gh)^{-1}} ] ~: g\in\Sigma^1, h\in\Sigma^1\cup\{1\} \}.$ Then $T$ is the direct sum of the  ideals
given in Theorem \ref{main1}-(1),
$$
T=\bigoplus_{[g]\in\Sigma^1/\sim}I_{[g]},
$$
being $[I_{[g]}, T, I_{[h]}]=[I_{[g]}, I_{[h]}, T]=[T, I_{[g]}, I_{[h]}]=0,$ if $[g]\neq [h].$
\end{cor}
\noindent {\bf Proof.} It follows as in \cite{CC2} Corollary 4.11.\qed

\section{The  simplicity and primeness of  $3-$Leibniz algebras} \setcounter{equation}{0}\ 

In this section we focus on the gr-simplicity of graded $3-$Leibniz algebras by centering
our attention in those of maximal length. This terminology is in a similar way than for graded Leibniz algebras \cite{CS1}. See also  \cite{C2, C3, CC1, CC2, CF, CO}.

\begin{DEF}\label{mal} We say that a graded $3-$Leibniz algebra $T=T_1\oplus(\bigoplus_{g\in\Sigma^1}T_g)$ is of {\em maximal length} if  $\dim T_g=1,$ for any $g\in\Sigma^1.$ 
\end{DEF}

\begin{lem}\label{4.2} Let $T=T_1\oplus(\bigoplus_{g\in\Sigma^1}T_g)$ be a graded $3-$Leibniz alebra of maximal length.  Then every ideal $I$ of $T$ decomposes as $I=(I\cap T_1)\oplus(\bigoplus_{g\in\Sigma^1}(I\cap T_g).$
\end{lem}
\noindent {\bf Proof.} See Lemma 5.2 in \cite{CO}.\qed

 From now on $T=T_1\oplus(\bigoplus_{g\in\Sigma^1}T_g)$ denotes a graded $3-$Leibniz alebra of maximal length. Using the previous lemma, we assert that given any nonzero ideal $I$ of $T$ then

\begin{eqnarray}\label{339}
I=(I\cap T_1)\oplus(\bigoplus_{g\in\Sigma^1_I} T_g),
\end{eqnarray}
where $\Sigma^1_I :=\{h\in\Sigma^1 : I\cap T_h\neq 0\}.$
In particular, case $I=J,$ we have
\begin{equation}\label{39'}
J=(J\cap T_1)\oplus(\bigoplus_{g\in\Sigma^1_J} T_g).
\end{equation}
From here, we can write
\begin{equation}\label{39''}
\Sigma^1=\Sigma^1_J\cup\Sigma^1_{J_0},
\end{equation}
where
$$
\Sigma^1_J :=\{g\in\Sigma^1 :T_g\cap J\neq 0\}=\{g\in\Sigma^1 :T_g\subset J\}, 
$$
and 
$$
\Sigma^1_{J_0} :=\{g\in\Sigma^1 : T_g\cap J=0\}.
$$
Therefore,
\begin{equation}\label{40}
T= T_1\oplus(\bigoplus_{h\in\Sigma^1_{J_0}}T_h)\oplus(\bigoplus_{g\in\Sigma^1_J}T_g).
\end{equation}
%The decomposition given by equation (\ref{40}) shows
%\begin{equation}\label{40'}
%T_1=\span_{\bbbf}\{[T_g, T_h, T_{(gh)^{-1}} ] ~: g\in\Sigma^1_{J_0}, %h\in\Sigma^1_{J_0}\cup\{1\} \}.
%\end{equation}

\begin{rem} The concept of connectivity of  elements in $\Sigma^1$ given in Definition \ref{conn} is not strong enough to detect if a given $g\in\Sigma^1$ belongs to $\Sigma^1_J$ or to $\Sigma^1_{J_0}.$  Consequently we lose the information respect to whether a given  space $T_g$ is contained in $J$ or not, which is fundamental to study the simplicity of $T.$ So, we are going to refine the concept of connections of  elements of $\Sigma^1$ in the set up of maximal length for graded $3-$Leibniz algebras. 
\end{rem}

\begin{DEF}\label{conn}
Let $g, ~h\in\Sigma^1_\Psi$ with $\Psi\in\{J, J_{0}\}.$ We say that {\em $g$ is $J-$connected to $h$} and denoted by $g\sim_{J}h$ if there exists a
family $\{g_1, g_2, g_3, ..., g_{2n+1}\}$  contained in $\pm\Sigma^1_{J_0}\cup\{1\}$ if $\Psi=J_{0}$ or contained in $\pm\Sigma^1_J\cup\{1\}$ if $\Psi=J,$ satisfying the following conditions;\\
\begin{itemize}
\item[(1)] $\{g_1, g_1g_2g_3, ..., g_1g_2g_3 ... g_{2n}g_{2n+1}\}\subset\pm\Sigma^1_\Psi,$\\
	
\item[(2)] $\{g_1g_2, g_1g_2g_3g_4, ... , g_1g_2g_3 ... g_{2n}\}\subset\pm\Sigma^0,$\\
	
\item[(3)] $g_1 = g$ and $g_1g_2g_3 ... g_{2n}g_{2n+1} \in\{h, h^{-1}\}.$
\end{itemize}
We will also that  $\{g_1, g_2, g_3, ..., g_{2n+1}\}$ is  a
$J-$connection  from $g$ to $h.$
\end{DEF}

The next result can be proved in a similar way to Proposition \ref{rel}.

\begin{pro}\label{rel2} The relation $\sim_J$ is an equivalence relation in both $\Sigma^1_{J}$ and $\Sigma^1_{J_0}.$
\end{pro}

Let us introduce the concepts of $\Sigma^1-$multiplicativity  in the framework  of  graded $3-$Leibniz algebra of maximal length, in a similar way to the ones for  graded Leibniz algebras and split $3-$Leibniz algebra in \cite{CD, CO}.

\begin{DEF}\label{malti} We say that a  graded $3-$Leibniz algebra of maximal length $T$ is {\em $\Sigma^1-$multiplicative} if 
\begin{itemize}
\item[(1)]Given $g\in\Sigma^1_{J_0}$ and $h, k\in\Sigma^1_{J_0}\cup\{1\}$ such that $gh\in\Sigma^0,~ghk\in\Sigma^1$ then $[T_g, T_h, T_k]\neq 0.$\\
	
\item[(2)] Given $g\in\Sigma^1_J$ and $h, k\in\Sigma^1\cup\{1\}$ such that $gh\in\Sigma^0,~ghk\in\Sigma^1_J$ then $[T_g, T_h, T_k]\neq 0.$
\end{itemize}
\end{DEF}

Another important notion related to   graded $3-$Leibniz algebra of maximal length $T$ is those of Lie-annihilator of $T= T_1\oplus(\bigoplus_{g\in\Sigma^1_J}T_g) \oplus(\bigoplus_{h\in\Sigma^1_{J_0}}T_h)$ (see Eq. (\ref{40}))  which define as follows;

\begin{DEF} The  Lie-annihilator of a   graded $3-$Leibniz algebra of maximal length $T$ is the set
\begin{eqnarray*}
Ann_{Lie}(T):&=&\{x\in T: [x,  T_1\oplus(\bigoplus_{h\in\Sigma^1_{J_0}}T_h), T_1\oplus(\bigoplus_{h\in\Sigma^1_{J_0}}T_h)]\\
&+&[T_1\oplus(\bigoplus_{h\in\Sigma^1_{J_0}}T_h), x,T_1\oplus(\bigoplus_{h\in\Sigma^1_{J_0}}T_h)]\\
&+&[T_1\oplus(\bigoplus_{h\in\Sigma^1_{J_0}}T_h), T_1\oplus(\bigoplus_{h\in\Sigma^1_{J_0}}T_h), x]=0 \}.
\end{eqnarray*}
\end{DEF}
Observe that $Ann(T)\subset Ann_{Lie}(T).$  Moreover, this is the natural ternary extension of the Lie annihilator of  graded
Leibniz algebras of maximal length, since in this case $\Sigma^1_J=0$ (See \cite{CD}). 

The symmetry of $\Sigma^1_J$ and $\Sigma^1_{J_0}$ will be understood as usual. That is, $\Sigma^1_\Psi,~~\Psi\in\{J, J_{0}\},$ is called symmetric if $g\in\Sigma^1_\Psi$ implies $g^{-1}\in\Sigma^1_\Psi.$  The same concept applies to the set $\Sigma^0.$ From now on we will assume that $\Sigma^1_{J_0}$ and $\Sigma^0$ are both symmetric.

\begin{pro}\label{lem1} Let $T$ be a graded $3-$Leibniz algebra of maximal length. Assume that $T$ is $\Sigma^1-$multiplicative, $Ann_{Lie}(T)=0$ and $T_1=\span_{\bbbf}\{[T_g, T_h, T_{(gh)^{-1}} ] ~: g\in\Sigma^1, h\in\Sigma^1\cup\{1\}.$ If $\Sigma^1_{J_0}$ has all of its elements $J-$connected, then any nonzero graded ideal $I$ of $T$ such that
$I\nsubseteq T_1+J$  satisfies that $I=T.$ 
\end{pro}
\noindent {\bf Proof.} Let $I$ be an ideal of $T$ such that $I\nsubseteq T_1+J.$ By Lemma \ref{4.2} and Eq. (\ref{40}), we can write
\begin{equation}\label{320}
I= (I\cap T_1)\oplus(\bigoplus_{h\in\Sigma^1_{I, {J_0}}} T_h)\oplus(\bigoplus_{g\in\Sigma^1_{I, J}}T_g).
\end{equation}
where $\Sigma^1_{I, {J_0}} :=\Sigma^1_I\cap\Sigma^1_{J_0}$ and $\Sigma^1_{I, J} :=\Sigma^1_I\cap\Sigma^1_{J}.$ Since $I\nsubseteq T_1+ J,$ from Eq. (\ref{320}) 
we have $I\cap(\bigoplus_{h\in\Sigma^1_{I, J_0}}T_h)\neq \{0\}.$ One gets $\Sigma^1_{I, J_0}\neq\emptyset$  and so we can take some $g_0\in\Sigma^1_{I, J_0}$ such that 
\begin{equation}\label{501}
0\neq T_{g_0}\subset I.
\end{equation}
Since $g_0$ is $J-$connected to any $h\in\Sigma^1_{J-0},$  we have a $J-$connection  $$\{g_1, g_2,...,
g_{2n+1}\}\subset\pm\Sigma^1_{J_0}\cup\{1\},$$ from $g_0$ to
$h$ satisfying the following conditions;\\
\begin{itemize}
\item[(1)] $\{g_1, g_1g_2g_3, ..., g_1g_2g_3 ... g_{2n}g_{2n+1}\}\subset\Sigma^1_{J_0},$\\
	
\item[(2)] $\{g_1g_2, g_1g_2g_3g_4, ... , g_1g_2g_3 ... g_{2n}\}\subset\Sigma^0,$\\
	
\item[(3)] $g_1 = g_0$ and $g_1g_2g_3 ... g_{2n}g_{2n+1} \in\{h, h^{-1}\}.$
\end{itemize}
Since $\Sigma^1_{J_0}$ and $\Sigma^0$ are both symmetric, we have  $g_1 g_2\in\Sigma^0$ and $g_1g_2g_3\in\Sigma^1_{j_0}.$ Then,  the $\Sigma^1-$multiplicativity and maximal length of $T$ allow us to get
$$
0\neq[T_{g_1}, T_{g_2}, T_{g_3}]=T_{g_1g_2g_3},
$$
so by Eq. (\ref{501}) and $T_{g_1}\subset I,$  we have
$$
0\neq T_{g_1g_2g_3}\subset I.
$$
We can argue in a similar way from $g_1 g_2 g_3\in\Sigma^0$ and $g_1g_2g_3g_4\in\Sigma^1_{j_0},$ to get
$$
0\neq T_{g_1g_2g_3g_4}\subset I.
$$
Following  this process with the $J-$connection $\{g_1, g_2,...,g_{2n+1}\},$ we obtain that
$$
0\neq T_{g_1g_2g_3 ... g_{2n}g_{2n+1}}\subset I, 
$$
and so  
\begin{equation}\label{313}
\hbox{either~~}~~0\neq T_{h}\subset I~~\hbox{ or~~}~~~
0\neq T_{h^{-1}}\subset I,~~~~~\forall h\in\Sigma^1_{J_0}.
\end{equation}

We claim that $J\cap T_1\subset Ann_{Lie}(T).$ Indeed, taking into account Eq. (\ref{101}) we have
$$
[T_1\oplus(\bigoplus_{h\in\Sigma^1_{J_0}}T_h), J\cap T_1, T_1\oplus(\bigoplus_{h\in\Sigma^1_{J_0}}T_h)]=
[T_1\oplus(\bigoplus_{h\in\Sigma^1_{J_0}}T_h), T_1\oplus(\bigoplus_{h\in\Sigma^1_{J_0}}T_h), J\cap T_1]=0.
$$
Now, for any $h\in\Sigma^1_{J_0},$ we have  $[J\cap T_1, T_1, T_h]\subset T_h\subset J,$ thanks to $J$ is an ideal. In the other hand, if $0\neq[J\cap T_1, T_1, T_h]$ then $h\in\Sigma^1_J,$ a contradiction. From here $[J\cap T_1, T_1, T_h]=0$ and so
\begin{equation}\label{315}
[J\cap T_1, T_1, \bigoplus_{h\in\Sigma^1_{J_0}}T_h]=0.
\end{equation}
Similarly $[J\cap T_1, \bigoplus_{h\in\Sigma^1_{J_0}}T_h, T_1].$ One can easily get $[J\cap T_1, T_1, T_1]=0.$ Therefore,
$$
[J\cap T_1,  T_1\oplus(\bigoplus_{h\in\Sigma^1_{J_0}}T_h), T_1\oplus(\bigoplus_{h\in\Sigma^1_{J_0}}T_h)]=0. 
$$
Hence 
\begin{equation}\label{314}
J\cap T_1\subset Ann_{Lie}(T)=0.
\end{equation}

Our next aim is to prove that
$$
T_1\subset I.
$$
Notice that Eq. (\ref{314}) produces  $T_1=\span_{\bbbf}\{[T_g, T_h, T_{(gh)^{-1}} ] ~: g\in\Sigma^1_{J_0}, h\in\Sigma^1_{J_0}\cup\{1\}.$ Let us study the products $[T_g, T_h, T_{(gh)^{-1}} ] ~: g\in\Sigma^1_{J_0}, h\in\Sigma^1_{J_0}\cup\{1\}$ in order to show $T_1\subset I.$ Suppose that $h=1$ then Eq. (\ref{313}) gives us $[T_g, T_1, T_{g^{-1}}]\subset I.$ Then we have that $g, h\in\Sigma^1_{J_0}.$ By $[T_g, T_h, T_{(gh)^{-1}} ]\neq 0$ and $gh\neq 1$ we get $gh\in\Sigma^0.$ Taking into account Eq. (\ref{313}) and $\Sigma^1-$multiplicative of $T,$ we have $[T_g, T_h, T_{(gh)^{-1}} ]\subset I.$ 

Finally, for any $h\in\Sigma^1,$ if $0\neq[T_g, T_h, T_{(gh)^{-1}} ]$ with $h\notin\{1, g^{-1}\}$ then by Eq. (\ref{313}), $\Sigma^1-$multiplicativity and maximal length, we get
$ 0\neq T_h=[T_g, T_h, T_{g^{-1}}]\subset I.$
\qed

\begin{pro}\label{lem3}  Let $T$ be a graded $3-$Leibniz algebra of maximal length. Assume that $T$ is $\Sigma^1-$multiplicative, $Ann_{Lie}(T)=0$ and $T_1=\span_{\bbbf}\{[T_g, T_h, T_{(gh)^{-1}} ] ~: g\in\Sigma^1, h\in\Sigma^1\cup\{1\}.$ If $\Sigma^1_J$ has all of its elements $J-$connected, then any nonzero graded ideal $I$ of $T$ such that $I\subset J$  satisfies that $I=J.$ 
\end{pro}
\noindent {\bf Proof.} By Eq. (\ref{40}) we can write
\begin{equation}\label{**}
	I= (I\cap T_1)\oplus(\bigoplus_{h_{j}\in\Sigma^1_{I, J}}T_{h_j})),
\end{equation}
where  $\Sigma^1_{I, J} :=\{g\in\Sigma^1_J~ :~T_g\cap I\neq 0\}=\{g\in\Sigma^1_J~ :~T_g\subset I\}.$ Note that, Eq. (\ref{314}) gives us
\begin{equation}\label{337}
J =\bigoplus_{h\in\Sigma^1_J}T_h.
\end{equation}
On the other hand $I \cap T_1\subset J\cap T_1 = 0,$ then we can rewrite Eq. (\ref{**}) by
$$
I =\bigoplus_{h_{j}\in\Sigma^1_{J, I}}T_{h_j}.
$$
Thus, we can find some  $h_{J, 0}\in\Sigma^1_{I, J}$ such that $0\neq T_{h_{J, 0}}\subset I.$

Observe that if $h\in\Sigma^1_J$ then $h^{-1}\notin\Sigma^1.$ Indeed, in the opposite case, the fact  $Ann_{Lie}(T)=0$ gives us  $[T_h, T_{h^{-1}}, T_g]\neq 0$ for some $g\in\Sigma^1_{J_0}.$  But then $g\in\Sigma^1_J,$ , a contradiction. Thus  $h^{-1}\notin\Sigma^1.$  Next, if $k\in\Sigma^1_J,$ a similar argument to the
one used in Proposition \ref{lem1} guarantees the existence of a $J-$connection $\{k_1, k_2,..., k_{2n}, k_{2n+1}\}\subset\pm\Sigma^1_J\cup\{1\}$ from $h_{J, 0}$ to $k$ such that 
$$
0\neq [[...[T_{h_{J, 0}}, T_{k_2}, T_{k_3}], ...]T_{k_{2n}}, T_{k_{2n+1}}]=T_k\subset I.
$$
From here and Eq. (\ref{337}) we have $J\subset I$.\qed

\begin{pro}\label{lem2}  Let $T$ be a graded $3-$Leibniz algebra of maximal length. Assume that $T$ is $\Sigma^1-$multiplicative, $Ann_{Lie}(T)=0$ and $T_1=\span_{\bbbf}\{[T_g, T_h, T_{(gh)^{-1}} ] ~: g\in\Sigma^1, h\in\Sigma^1\cup\{1\}.$ If $\Sigma^1_J$ has all of its elements $J-$connected, then any nonzero graded ideal $I$ of $T$ such that
$I\subset J$  satisfies that $I=J$ or $J=I\oplus K$ with $K$ a graded
ideal of $T.$ 
\end{pro}
\noindent {\bf Proof.} By Eqs. (\ref{39''}),  (\ref{40}) and $I\subset J$ we can write
\begin{equation}\label{*}
I= (I\cap T_1)\oplus(\bigoplus_{h_{j}\in\Sigma^1_{I, J}}T_{h_j}),
\end{equation}
where  $\Sigma^1_{I, J} :=\Sigma^1_I\cap\Sigma^1_J.$ Eq. (\ref{314}) gives us
 $I \cap T_1\subset J\cap T_1 = 0,$  then we can rewrite Eq. (\ref{*}) by $I =\bigoplus_{h_{j}\in\Sigma^1_{I, J}}T_{h_j},$ with $\Sigma^1_{I, J}\neq\emptyset.$ Thus, we can take $h_0\in\Sigma^1_{I, J}$ such that
$$
0\neq T_{h_0}\subset I.
$$
We can argue with the $\Sigma^1-$multiplicativity and the maximal length of $T$ as in Proposition \ref{lem1} to conclude that given any $h\in\Sigma^1_J,$ there exists a $J-$connection $\{h_0=g_1, g_2, ... , g_{2n+1}\}\subset \pm\Sigma^1_{J_0}\cup\{1\}$ from $h_0$ to $h$ such that
\begin{equation*}
0\neq [[...[T_{h_0}, T_{g_2}, T_{g_3}], ...], T_{g_{2n}}, T_{g_{2n+1}}] = T_k\subset I,
\end{equation*}
for some $k\in\{h, h^{-1}\}$
That is, $k\in\Sigma^1_{I, J}$ for some $k\in\{h, h^{-1}\}$ and for all $h\in\Sigma^1_J.$
Suppose $h^{-1}_0\in\Sigma^1_{I, J}.$ We get that $\{{h_0}^{-1}={g_1}^{-1},  {g_2}^{-1}, ..., {g_{2n+1}}^{-1}\}\subset \pm\Sigma^1_{J_0}\cup\{1\}$ is a $J-$connection from ${h_0}^{-1}$ to $h$ such that
\begin{equation}\label{316}
0\neq [[...[T_{{h_0}^{-1}}, T_{{ g_2}^{-1}}, T_{{g_3}^{-1}}], ...], T_{{g_{2n}}^{-1}}, T_{{g_{2n+1}}^{-1}}] = T_{k^{-1}}\subset I,
\end{equation}
and so $\LL_h+\LL_{h^{-1}}\subset I$ for any $h\in\Sigma^1_J.$  Thus,  Eqs.  (\ref{337}) imply that $J\subset I.$ So $I=J$

Next, suppose there is not any $h_0\in\Sigma^1_{I, J}$ such that $h_0^{-1}\in\Sigma^1_{I,J}.$ Taking account into Eq. (\ref{316}) and the fact that $I\subset J$ allows us to
write $\Sigma^1_J=\Sigma^1_{I, J}\cup(\Sigma^1_{I, J})^{-1},$ where $(\Sigma^1_{I, J})^{-1}~:=\{h^{-1}\in\Sigma^1_J ~:~h\in\Sigma^1_{I, J}\}.$
By Eqs. (\ref{40}) and (\ref{314}),  denote $K=\bigoplus_{h\in\Sigma^1_{I, J}}T_{h^{-1}}),$ and we get
\begin{equation}\label{55'}
J=I\oplus K.
\end{equation}
It remains to show that $K$ is a graded ideal of $T.$  Note that, by Eq. (\ref{101}) we have $[T, T, K]=[T, K, T]=0,$ and also 
\begin{eqnarray}\label{56}
\nonumber[K, T, T]&=& [K, T_1\oplus(\bigoplus_{g\in\Sigma^1_J}T_g) \oplus(\bigoplus_{h\in\Sigma^1_{J_0}}T_h), T_1\oplus(\bigoplus_{g\in\Sigma^1_J}T_g) \oplus(\bigoplus_{h\in\Sigma^1_{J_0}}T_h)]\\
&\subset&[K, T_1, T_1]+[K,T_1, \bigoplus_{g\in\Sigma^1_J}T_g]+[K, T_1, \bigoplus_{h\in\Sigma^1_{J_0}}T_h]\\
\nonumber&+&[K, \bigoplus_{g\in\Sigma^1_J}T_g, T_1]+[K, \bigoplus_{g\in\Sigma^1_J}T_g, \bigoplus_{g'\in\Sigma^1_J}T_{g'}]+[K, \bigoplus_{g\in\Sigma^1_J}T_g, \bigoplus_{h\in\Sigma^1_{J_0}}T_h]\\
\nonumber&+&[K, \bigoplus_{h\in\Sigma^1_{J_0}}T_h, T_1]+[K, \bigoplus_{h\in\Sigma^1_{J_0}}T_h, \bigoplus_{g\in\Sigma^1_J}T_g]+[K, \bigoplus_{h\in\Sigma^1_{J_0}}T_h, \bigoplus_{h'\in\Sigma^1_{J_0}}T_{h'}].	
\end{eqnarray}
It is clear that
\begin{equation*}
[K, T_1, \bigoplus_{g\in\Sigma^1_J}T_g]=[K, \bigoplus_{h\in\Sigma^1_{J_0}}T_h, \bigoplus_{g\in\Sigma^1_J}T_g]=[K, \bigoplus_{g\in\Sigma^1_J}T_g, T_1]=0,
\end{equation*}
and also
$$
[K, \bigoplus_{g\in\Sigma^1_J}T_g, \bigoplus_{h\in\Sigma^1_{J_0}}T_h] =0,~~~[K, T_1, T_1]\subset K.
$$
Let us consider the  summand $[K, T_1, \bigoplus_{h\in\Sigma^1_{J_0}}T_h]$ in Eq. (\ref{56}), and suppose there exist $g\in\Sigma^1_{I, J}$ and $h\in\Sigma^1_{J_0}$ such that $0\neq[T_{g^{-1}}, T_1, T_h].$ Since $T_{g^{-1}}\subset K\subset J,$ we get $g^{-1} h\in\Sigma^1_J.$ By the $\Sigma^1-$multiplicativity of $T,$ the symmetries of $\Sigma^1_J,~~\Sigma^1_{J_0},$ and the fact $T_g\subset I,$ one gets $0\neq[T_{g^{-1}}, T_1, T_h]=T_{g^{-1}h}\subset I.$ That is, $g^{-1}h\in\Sigma^1_{I, J}.$ Hence, $gh^{-1}\in(\Sigma^1_{I, J})^{-1},$ and so $[T_g, T_1, T_h^{-1}]=T_{gh^{-1}}\subset K.$  Similarly, we also get $[K, \bigoplus_{h\in\Sigma^1_{J_0}}T_h, T_1]\subset K.$

Finally, consider the summand $[K, \bigoplus_{h\in\Sigma^1_{J_0}}T_h, \bigoplus_{h'\in\Sigma^1_{J_0}}T_{h'}]$ in Eq. (\ref{56}),  and suppose there exist $g\in\Sigma^1_{I, J}$ and $h, h'\in\Sigma^1_{J_0}$ such that $0\neq[T_{g^{-1}}, T_h, T_{h'}].$ Since $T_{g^{-1}}\subset K\subset J,$ we get $g^{-1} h h'\in\Sigma^1_J.$ By the $\Sigma^1-$multiplicativity of $T,$ the symmetries of $\Sigma^1_J,~~\Sigma^1_{J_0},$ and the fact $T_g\subset I,$ one gets  $0\neq[T_{g^{-1}}, T_h, T_{h'}]=T_{g^{-1}hh'}\subset I.$ That is, $g^{-1}hh'\in\Sigma^1_{J, I}.$ Hence, $gh^{-1}h'^{-1}\in(\Sigma^1_{J, I})^{-1},$ and so $[T_g, T_h^{-1}, T_{h'^{-1}}]=T_{gh^{-1}h'^{-1}}\subset K.$ Consequently, $K$ is
a graded ideal of $T.$\qed

We introduce the definition of gr-primeness in the framework of graded $3-$Leibniz
algebras following the same motivation in the case of gr-primeness of graded Liebniz algebras (see \cite{CD}).

\begin{DEF}\label{prim} A graded $3-$Leibniz algebra $T$ is said to be {\em gr-prime} if given two ideals $I$ and $K$ of $T$ satisfying $[I, K, I] + [K, I, I] + [I, I, K] = 0,$ then either $I\in\{0, J, T\}$ or $K \in\{0, J, T\}.$
\end{DEF}

We also note that the above Definition agrees with the Definition of prime graded $3-$Lie algebra, since $J = 0$ in this case.

\begin{cor}\label{4.12} Under the hypotheses of Proposition \ref{lem2}, if furthermore $T$ is gr-prime, then any non-zero graded ideal $I$ of $T$ such that $I\subseteq J$ satisfies $I=J.$
\end{cor}
\noindent {\bf Proof.} By Proposition \ref{lem2}, we have $J=I\oplus K$ with $I,~~K$ graded ideals of $T,$ being
$$
[I, K, I]+[K, I, I]+[I, I, K]=0,.
$$
as consequence of $I,~~K\subset J.$ Now, the primeness of $T$ completes the proof.\qed

In the following, consider $T$ is a graded $3-$Leiniz algebra with $Ann(T)=0$ and  $T_1=\span_{\bbbf}\{[T_g, T_h, T_{(gh)^{-1}} ] ~: g\in\Sigma^1, h\in\Sigma^1\cup\{1\} \}.$ For any $g\in\Sigma^1_\Psi,~~\Psi\in\{J, J_0\},$ we denote by  
$$
\Lam^g_\Psi :=\{h\in\Sigma^1_\Psi : h\sim_J g\}.
$$
Our goal is to introduce an ideal of $T$ associated to $\Lam^g_\Psi.$ For a fixed  $g\in\Sigma^1_\Psi,~~\Psi\in\{J, J_0\},$ we define a  linear space $I_{1, \Lam^g_\Psi}$ by 
\begin{equation}\label{119}
I_{1, \Lam^g_\Psi}= \span_{\bbbf}\{[T_h, T_k, T_{(hk)^{-1}} ] ~: h\in\Lam^g_\Psi, k\in \Lam^g_\Psi\cup\{1\} \}\sub T_1.
\end{equation}
Next, we define
\begin{equation}\label{119'}
\vv_{\Lam^g_\Psi}=\bigoplus_{h\in\Lam^g_\Psi} T_h.
\end{equation}
We denote by $I_{\Lam^g_\Psi}$ the following subspace of $T,$ 
\begin{equation}\label{119'}
I_{\Lam^g_\Psi}=I_{1, \Lam^g_\Psi}\oplus\vv_{\Lam^g_\Psi}.
\end{equation}

\begin{pro}\label{4.14}  If   $Ann(T)=0$ and  $T_1=\span_{\bbbf}\{[T_g, T_h, T_{(gh)^{-1}} ] ~: g\in\Sigma^1, h\in\Sigma^1\cup\{1\} \}.$ Then $I_{\Lam^g_J}$ is a graded ideal of $T$  for any $g\in\Sigma^1_J.$
\end{pro}
\noindent {\bf Proof.} From 
\begin{eqnarray*}
T_1&=&\span_{\bbbf}\{[T_g, T_h, T_{(gh)^{-1}} ] ~: g\in\Sigma^1, h\in\Sigma^1\cup\{1\} \}\\
&=&\span_{\bbbf}\{[T_g, T_h, T_{(gh)^{-1}} ] ~: g\in\Sigma^1_{J_0}, h\in\Sigma^1_{J_0}\cup\{1\} \},
\end{eqnarray*}
we get $I_{1, \Lam^g_J}=0$ and so
\begin{equation}\label{66'}
I_{\Lam^g_J}=\vv_{\Lam^g_J}=\bigoplus_{h\in\Lam^g_J} T_h.
\end{equation}
We are going to show that $[T, I_{\Lam^g_J}, T]+[T, T, I_{\Lam^g_J}]+[I_{\Lam^g_J}, T, T]\subset I_{\Lam^g_J}.$ Observe that, by Eq. (\ref{101}) we have 
\begin{equation*}
[T, I_{\Lam^g_J}, T]=[T, T, I_{\Lam^g_J}]=0.
\end{equation*}
It remains to verify that 
$$
[I_{\Lam^g_J}, T, T]\subset I_{\Lam^g_J},~~~\hbox{~~for~~any~}g\in\Sigma^g_J.
$$
Using Eq. (\ref{40}), we have
\begin{eqnarray}\label{325}
\nonumber[I_{\Lam^g_J}, T, T]&=& [I_{\Lam^g_J}, T_1\oplus(\bigoplus_{g\in\Sigma^1_J}T_g) \oplus(\bigoplus_{h\in\Sigma^1_{J_0}}T_h), T_1\oplus(\bigoplus_{g\in\Sigma^1_J}T_g) \oplus(\bigoplus_{h\in\Sigma^1_{J_0}}T_h)]\\
&\subset&[I_{\Lam^g_J}, T_1, T_1]+[I_{\Lam^g_J}, T_1, \bigoplus_{g\in\Sigma^1_J}T_g]
+[I_{\Lam^g_J}, T_1, \bigoplus_{h\in\Sigma^1_{J_0}}T_h]\\
\nonumber&+&[I_{\Lam^g_J}, \bigoplus_{g\in\Sigma^1_J}T_g, T_1]
+[I_{\Lam^g_J}, \bigoplus_{g\in\Sigma^1_J}T_g, \bigoplus_{g'\in\Sigma^1_J}T_{g'}]\\
\nonumber&+&[I_{\Lam^g_J}, \bigoplus_{g\in\Sigma^1_J}T_g, \bigoplus_{h\in\Sigma^1_{J_0}}T_h]
+[I_{\Lam^g_J}, \bigoplus_{h\in\Sigma^1_{J_0}}T_h, T_1]\\
\nonumber&+&[I_{\Lam^g_J}, \bigoplus_{h\in\Sigma^1_{J_0}}T_h, \bigoplus_{g\in\Sigma^1_J}T_g]+[I_{\Lam^g_J}, \bigoplus_{h\in\Sigma^1_{J_0}}T_h, \bigoplus_{h'\in\Sigma^1_{J_0}}T_{h'}].	
\end{eqnarray}
Here, it is clear that
\begin{equation*}
[I_{\Lam^g_J}, T_1, \bigoplus_{g\in\Sigma^1_J}T_g]=[I_{\Lam^g_J}, \bigoplus_{h\in\Sigma^1_{J_0}}T_h, \bigoplus_{g\in\Sigma^1_J}T_g]=[I_{\Lam^g_J}, \bigoplus_{g\in\Sigma^1_J}T_g, T_1]=0,
\end{equation*}
and also
$$
[I_{\Lam^g_J}, \bigoplus_{g\in\Sigma^1_J}T_g, \bigoplus_{h\in\Sigma^1_{J_0}}T_h] =0,~~~[I_{\Lam^g_J}, T_1, T_1]\subset I_{\Lam^g_J}.
$$
Let us consider the  summand $[I_{\Lam^g_J},  \bigoplus_{h\in\Sigma^1_{J_0}}T_h, T_1]$ in Eq. (\ref{325}). Given any $k\in\Lam^g_J$ and $h\in\Sigma^1_{J_0}$ such that $[T_k, T_h, T_1]\neq 0.$ If $kh\neq 1,$ then we have $kh\in\Sigma^1_J$ and so $\{k, h, 1\}$ is a $J-$connection from $k$ to $kh.$ By symmetry and transitivity of $\sim_J$ in $\Sigma^1_{J_0},$ we get $kh\in\Lam^g_J.$ Hence $[I_{\Lam^g_J},  \bigoplus_{h\in\Sigma^1_{J_0}}T_h, T_1]\subset  I_{\Lam^g_J}.$ Similarly, we also get $[I_{\Lam^g_J}, T_1,  \bigoplus_{h\in\Sigma^1_{J_0}}T_h]\subset  I_{\Lam^g_J}.$

Finally, consider the summand $[I_{\Lam^g_J}, \bigoplus_{h\in\Sigma^1_{J_0}}T_h, \bigoplus_{h'\in\Sigma^1_{J_0}}T_{h'}]$ in Eq. (\ref{325}),  and suppose there exist $k\in\Lam^g_J$ and $h, h'\in\Sigma^1_{J_0}$ such that $0\neq[T_k, T_h, T_{h'}].$ By $khh'\neq 1,$ we have $khh'\in\Sigma^1_J$ and so $\{k, h, h'\}$ is a $J-$connection from $k$ to $khh'.$ By symmetry and transitivity of $\sim_J$ in $\Sigma^1_{J_0},$ we get $khh'\in\Lam^g_J.$ Hence $[I_{\Lam^g_J}, \bigoplus_{h\in\Sigma^1_{J_0}}T_h, \bigoplus_{h'\in\Sigma^1_{J_0}}T_{h'}]\subset  I_{\Lam^g_J}.$

Consequently, $[T, I_{\Lam^g_J}, T]+[T, T, I_{\Lam^g_J}]+[I_{\Lam^g_J}, T, T]\subset I_{\Lam^g_J}.$ So $I_{\Lam^g_J}$ is a graded ideal of $T.$\qed

Now, we are ready to state our main result;

\begin{thm}\label{Main} Suppose $T$ is a graded $3-$Leibniz algebra of maximal length.  Assume that $T$ is $\Sigma^1-$multiplicative, $T_1=\span_{\bbbf}\{[T_g, T_h, T_{(gh)^{-1}} ] ~: g\in\Sigma^1, h\in\Sigma^1\cup\{1\}\}$ and $Ann_{Lie}(T)=0.$ Then $T$ is gr-simple if and only  if it is gr-prime and $\Sigma^1_J,~~\Sigma^1_{J_0},$ have all of their  elements $J-$connected.
\end{thm}
\noindent {\bf Proof.} Suppose $T$ is a gr-simple $3-$Leibniz algebra. If $\Sigma^1_J\neq\emptyset$ and we take $g\in\Sigma^1_J.$ Lemma \ref{4.14} gives us  $I_{\Lam^g_J}$ is a nonzero graded ideal of $T.$ By simplicity of $T,$ we have $I_{\Lam^g_J}=J$ and $Z_{Lie}(T)=0$ implies that $J\cap T_1=0,$ so
$$
I_{\Lam^g_J}=J=\bigoplus_{h_j\in{\Sigma^1_J}}T_{h_j}.
$$
Hence, $\Sigma^1_J=\Lam^g_J$ and $\Sigma^1_J$ has all of its elements $J-$connected. 

Consider now any $h\in\Sigma^1_{J_0}$ and the graded subspace $I_{\Lam^h_{J_0}}.$ Denote by $<I_{\Lam^h_{J_0}}>$ the ideal of $T$ generated by $I_{\Lam^h_{J_0}}.$  Since $J$ is a graded ideal of $T,$ we have $<I_{\Lam^h_{J_0}}>\cap\bigoplus_{k\in\Sigma^1_{J_0}}T_k$ is contained in the linear span of the set
\begin{itemize}
\item[~]	$\{[[...[x_{h'}, x_{g_1}, x_{g_2}]...], x_{g_{2n}}, x_{g_{2n+1}}]; [ x_{g_{2n+1}}, x_{g_{2n}}, [...[x_{g_2}, x_{g_2}, x_{h'}]...]];$\\
\item[~] $ [[...[ x_{g_1}, x_{g_2}, x_{h'}], ...], x_{g_{2n}}, x_{g_{2n+1}}]; [x_{g_{2n+1}}, x_{g_{2n}}, [...[x_{h'}, x_{g_2}, x_{g_1}], ...]];$\\
\item[~] $ [[...[x_{g_1}, x_{h'}, x_{g_2}], ...], x_{g_{2n}}, x_{g_{2n+1}}]; [x_{g_{2n+1}}, x_{g_{2n}}, [...[x_{g_2}, x_{g_1}, x_{h'}], ...]]\}$, 
\end{itemize}
where $0\neq x_{h'}\in I_{\Lam^h_{J_0}}, 0\neq x_{g_i}\in T_{g_i},~~g_i\in\Sigma^1_{J_0},$ and $n\in\bbbn.$ 

By simplicity of $T$ we get $<I_{\Lam^h_{J_0}}>=T.$ From here, given any $l\in\Sigma^1_{J_0}$ the above observation gives us that we can write $l=h'g_1g_2....g_{2n}g_{2n+1}$ for $h'\in\Lam^h_{J_0}, g_i\in\Sigma^1_{J_0}$ and being the partial sums non zero. Hence, $\{h', g_1, g_2, ..., g_{2n+1}\}$ is a $J-$connection from $h'$ to $l.$ By the symmetry and transitivity of $\sim_J$ in $\Sigma^1_{J_0},$ we deduce that $h$ is $J-$connection to any $l\in\Sigma^1_{J_0}.$ Consequently, Proposition \ref{rel2} gives us that $\Sigma^1_{J_0}$ has all of its elements $J-$connected. Finally, the gr-simplicity of $T$ implies that $T$ is gr-prime.

Conversely, given $I$ a nonzero ideal of $T,$  we distinguish two cases:\\

{\bf Case 1.} If $I\nsubseteq T_1+J,$ then Proposition \ref{lem1} gives us $I=T.$\\ 

{\bf Case 2.} If $I\subset T_1+J,$ then $I\cap T_1\subset Ann_{Lie}(T).$ In fact, for any $x\in I\cap T_1,$ and $g\in\Sigma^1_{J_0}$ we have $[x, T_1, T_g]\subset J$ which implies that $[x, T_1, T_g]=0.$ Thus $I\cap T_1\subset Ann_{Lie}(T)=0.$
Applying Lemma \ref{4.2}  we get that $I\subset J,$ and Proposition \ref{lem3} gives us  $I=J.$ 

Therefore, we have proved that either $I=T$ or $I=J,$ which shows that $T$ is gr-simle.\qed

\section{Example} \setcounter{equation}{0}\

In this section, we provide an example to clarify the results in this article. The process of the example is described in two steps as follows.\\ 

{\bf Step 1.} (Complex $5-$dimentional graed $3-$Leibniz algebra $T$)  Consider the complex vector space $T$ with the basis $\{e, h, f, p, q\}$ defined by the following multiplication:

\begin{equation*}
\begin{split}
&[e, h] = 2e,~~~~ [h, f] = 2f,~~~~ [e, f] = h,\\
&[h, e] = -2e,~~~~[f, h] =-2f,~~~~ [f, e]= -h,\\
&[p, h] = p,~~~~ [p, f] = q,\\
&[q, h] = -q,~~~~ [q, e] =-p,
\end{split}
\end{equation*}
where omitted products are equal to zero. Then $(T, [., .])$ is a (non Lie) Leibniz algebra (see Example 3.2 in \cite{Om}). Next, we constract a $3-$Leibniz algebra which introduced in \cite{CLP}. We define   the operation $[., ., .]~:~T\times T\times T\longrightarrow T$ given by $[x, y, z]~:=[x, [y, z]].$ Then $(T, [., ., .])$  is also a complex $3-$Leibniz algebra, thanks to Proposition 3.2 in \cite{CLP}. Moer precisely, the triple product on $T$ defined by
\begin{equation*}
\begin{split}
&[h, h, e] = 4e,~~~ [h, h, f] = 4f,~~~ [f, h, e]=2h,~~~~~[e, h, f]=2h,\\
&[e, e, f]=2e,~~~[f, f, e]=2f,~~~~~[h, e, h] = -4e,~~~ [h, f, h] = -4\\
&[f, e, h] =-2h,~~~~[e, f, h]=-2h,~~~~[e, f, e]=-2e,~~~~ [f, e, f]=-2f\\
&[q, h, e]=2p,~~~~  [p, e, f] = p,~~~~~ [q, f, e]=q,~~~~~~ [p, h, f]=2q\\
&[q, e, h] =-2p,~~~~~ [p, f, e]=-p,~~~~~~  [q, e, f]=-q,~~~~~  [p, f, h]=-2q,
\end{split}
\end{equation*}
	
where omitted products are equal to zero. The $3-$Leibniz algebra $T$ can be $\bbbz-$graded as $T=\bigoplus_{z\in\bbbz}T_z,$ where 
$$
T_0=<h>,~ T_1=<p>,~ T_{-1}=<q>,~ T_2=<e>,~ T_{-2}=<f>,
$$ 
and $ T_z=0$ for any $z\notin\{0, \pm1, \pm2\}.$ The standaed embedding algebra of $T$ is $2-$graded algebra $\aa=(0, T).$ By definition, the support of the grading is $\Sigma^1=\{-2, -1, 1, 2\}$ an also $\Sigma^0=\emptyset.$	

The graded ideals of $T$ are $\{0\}, T$ and $J=T_{-1}\oplus T_1.$ One can check that the graded ideal $J$  satisfies in Eq. (\ref{101}), and is the only nonzero improper graded ideal of $T.$ Then $T$ is a gr-simple $\bbbz-$graded $3-$Leibnize algebra. Note that $Ann(T)=0.$\\

{\bf Step 2.} (Clarify the results in sections 3 and 4) By step 1, we have
$$
T=T_0\oplus(\bigoplus_{z\in\Sigma^1}T_z),
$$  
where the support of the grading is $\Sigma^1=\{-2, -1, 1, 2\},$ which is symmetric. By Theorem \ref{main1}, since $T$ is gr-simple we must have  $\Sigma^1$ has all of its elements connected. In fact, by Definition \ref{conn} we have\\

\begin{itemize}
\item[-] $\{1, 2, -1\}$ is a connection from $1$ to $2.$ Indded, $g_1=1$ and $g_1g_2g_3=1+2+(-1)=2\in\Sigma^1, ~~g_1g_2g_3=h.$ 

\item[-]Similarly, $\{-1, -2, 1\}$ is a connection from $-1$ to $2.$
\end{itemize}
Then  we get that $[1]=\{\pm 1, \pm 2\}=\Sigma^1.$ The tabel of triple product gives that $T_1=<[T_g, T_h, T_{(gh)^{-1}}]>$ where  $g\in\Sigma^1, h\in\Sigma^1\cup\{1\}.$ The decomposition of $T$ is the form
$$
T=T_1\oplus(\bigoplus_{g\in\Sigma^1}T_g)=\uu\oplus\sum_{[g]\in\Sigma^1/\sim}I_{[g]}.
$$
Taking into account $[1]=\Sigma^1$ we get $T$ is the direct sum of the graded ideals
given in Theorem \ref{main1}-(1).

It is clear that $T$ is a $\bbbz-$graded $3-$Leibniz algebra of maximal length. By graded ideal $J=T_{-1}\oplus T_1,$ we have $\Sigma^1_J=\{-1, 1\}$ and $\Sigma^1_{J_0}=\{-2, 2\}.$ From here and the tabel of triple product, it is easy to see that $T$ is $\Sigma^1-$multiplicative. We also have $Ann_{Lie}(T)=0.$ 

Note that, $J$ is the only nonzero graded ideal of $T$  and the tabel of triple prouduct in $T$ shows that $T$ is a gr-prim $3-$Leibniz algebra (see Definition \ref{prim}). By $\Sigma^1=\Sigma^1_J\cup\Sigma^1_{J_0}$ we get that both $\Sigma^1_J,~~\Sigma^1_{J_0},$ have all of their  elements $\Sigma^1_{J_0}-$connected. From here, we conclude that $T$ is gr-simple if and only  if it is gr-prime and $\Sigma^1_J,~~\Sigma^1_{J_0},$ have all of their  elements $\Sigma^1_{J_0}-$connected. This verifies Teorem \ref{Main} in section 4.\\

Data availability statement; My manuscript has no associated data

\end{document}